\documentclass[a4paper,12pt]{article}
\usepackage{amssymb,amsmath}
\newcommand{\Q}{\mathbb{Q}}
\newcommand{\F}{\mathbb{F}}
\newcommand{\R}{\mathbb{R}}
\newcommand{\C}{\mathbb{C}}
\newcommand{\N}{\mathbb{N}}
\newcommand{\Proj}{\mathbb{P}}
\renewcommand{\b}[1]{{\bf #1}}
\newcommand{\x}{\b{x}}
\newcommand{\Z}{\mathbb{Z}}
\newtheorem{theorem}{Theorem}
\newtheorem{conjecture}{Conjecture}

\begin{document}
\title{$p$-adic Zeros of Systems of Quadratic Forms}
\author{D.R. Heath-Brown\\Mathematical Institute, Oxford}
\date{}
\maketitle

This survey concerns the following problem.  Let $K$ be a field and
let $r\in\N$.  Define $\beta(r;K)$ to be the largest integer $n$ for
which there exist quadratic forms 
\[q_i(x_1,\ldots,x_n)\in K[x_1,\ldots,x_n]\;\;\;(1\le i\le r)\]
having only the trivial common zero over $K$. Thus $\beta(r;K)$ is
also the smallest integer such that any such system in at least
$1+\beta(r;K)$ variables has a non-trivial common zero. When $r=1$ the
number $\beta(1;K)$ is the $u$-invariant of the field $K$.

As examples one has $\beta(1;\R)=\infty$, $\beta(1;\C)=1$ (since the
form $x_1^2$ has only the trivial zero $x_1=0$), and indeed
$\beta(r;\C)=r$.

We will be primarily interested in the case in which $K$ is a
$p$-adic field $\Q_p$.  It is well known in this case that
$\beta(1;\Q_p)=4$.  For example if $p$ is an odd prime and $k$ is a
quadratic non-residue of $p$ then $x_1^2-kx_2^2+p(x_3^2-kx_4^2)$ has
only the trivial zero over $\Q_p$, while any form in 5 variables has a
non-trivial zero.  Thus the key question is what one can say about
$\beta(r;\Q_p)$ in general.

Why should one be interested in such problems? Firstly, systems of
quadratics are fundamental to Diophantine analysis, since any
Diophantine equation may be reduced to such a system. Secondly, in
certain circumstances there are local-to-global principles for such
systems.  For example, suppose we have a system of quadratic forms
over $\Q$ which defines a smooth variety, and suppose also 
that the number of variables $n$ exceeds $2r^2+2r$. Then a
theorem of Birch \cite{birch}, proved via the Hardy--Littlewood circle
method, shows that there is a non-trivial rational point provided that
there is a non-trivial point over every completion of $\Q$.  It is
clear that one cannot drop the condition over $\R$, but it is natural
to ask whether the $p$-adic conditions are satisfied automatically for
$n\ge 2r^2+2r$. Thirdly,
the reduction of general Diophantine equations to systems of
quadratics can be made sufficiently efficient in special circumstances
that information about values of $\beta(r;\Q_p)$ can yield worthwhile
information about higher degree equations.  An example of this occurs
in the author's work \cite{hbquartics}, where it is shown that for any
prime $p$ different from 2 or 5, a
quartic form over $\Q_p$ in $n$ variables has a nontrivial zero,
provided that $n>16+\beta(8;\Q_p)$.  One would therefore like to
know how many variables are needed for a system of 8 quadratic forms
to have a non-trivial zero.

What might one expect about $\beta(r;\Q_p)$ ?  It was conjectured by
Artin, see \cite[p.x]{Artin}, that a $p$-adic form of degree $d$ in $n$
variables should have a non-trivial zero as soon as $n>d^2$.  A
consequence of this would be that a system of $r$ quadratic forms
in $n$ variables would have a non-trivial zero as soon as $n>4r$.  On
the other hand one may easily construct systems in $4r$ variables
having only the trivial zero.  Indeed if $q(x_1,x_2,x_3,x_4)$ has only
the trivial zero then one may take 
\[q_1=q(x_1,\ldots,x_4),\; q_2=q(x_5,\ldots,x_8),\;\ldots\;,
q_r=q(x_{4r-3},\ldots,x_{4r}).\]  
Thus one is led to the following conjecture.
\begin{conjecture}
For any $r\in\N$ and any prime number $p$ one has $\beta(r;\Q_p)=4r$.
\end{conjecture}

In fact Artin's Conjecture is known to be false (Terjanian
\cite{Ter}), but none of the known counterexamples relate to systems
of quadratic forms.  Thus the conjecture above remains open.

The most important result known on Artin's Conjecture is probably that
of Ax and Kochen \cite{AK}, who showed that for any given degree $d$
there is a corresponding $p(d)$ such that the conjecture is true for
primes $p\ge p(d)$.  One may deduce that for any $r$ there is a
corresponding $p'(r)$ such that $\beta(r;\Q_p)=4r$ as soon as $p\ge p'(r)$.

For small integers $r$ more has been proved.  As has already been remarked,
one has $\beta(1;\Q_p)=4$ for every prime $p$.  This was known
implicitly in the 19th century, but was proved explicitly by Hasse
\cite{Has}.  For $r=2$ one similarly has $\beta(2;\Q_p)=8$ for every
prime $p$, as was established by Demyanov \cite{Dem} in 1956.  However even
for $r=3$ the picture is incomplete.  Here it was shown by Schuur
\cite{Schuur}, building on work of Birch and Lewis \cite{BL}, that
$\beta(3;\Q_p)=12$ provided that $p\ge 11$.  Thus we have the
following problem.\bigskip

{\bf Open Question} {\em Is it true that} $\beta(3;\Q_p)=12$ {\em for
  all primes} $p$ ?\bigskip

There are two major lines of attack on such questions.  The first
traces its roots through work of Birch, Lewis and Murphy \cite{BLM}
(1962), Birch and Lewis \cite{BL} (1965), and Schmidt \cite{sch} (1980).
The basic idea is to choose representatives $q_i$ for the system 
in such a way that $q_i(x_1,\ldots,x_n)\in\Z_p[x_1,\ldots,x_n]$. Reduction
modulo $p$ then yields forms
$Q_i(x_1,\ldots,x_n)\in\F_p[x_1,\ldots,x_n]$, say.  Then, if
the system $Q_1,\ldots,Q_r$ has a non-singular zero over $\F_p$
one can lift it to a non-trivial zero of $q_1,\ldots,q_r$ over $\Z_p$,
by Hensel's Lemma.

By the Chevalley--Warning Theorem the system $Q_1,\ldots,Q_r$
certainly has a non-trivial zero over $\F_p$ as soon as $n>2r$, so the
key issue is whether or not we can produce a non-singular zero.  This
is clearly not possible in general.  Indeed nothing that has been said
so far precludes the possibility that the forms $Q_i$ all vanish.
Thus the strategy is to start by choosing a good integral model for
the system $<q_1,\ldots,q_r>_{\Q_p}$, by removing as many excess
factors of $p$ as one can. Here, when we refer to a good model for the
system, one should observe that one may make invertible linear changes of
variable 
\begin{equation}\label{xmap}
\x\mapsto M\x,\;\;\;\;\; M\in {\rm GL}_n(\Q_p)
\end{equation}
amongst $x_1,\ldots,x_n$, and invertible linear changes 
\begin{equation}\label{qmap}
\mathbf{q}\mapsto P\mathbf{q},\;\;\;\;\; P\in {\rm GL}_r(\Q_p)
\end{equation}
amongst the forms
$q_1,\dots,q_r$, without affecting the existence or otherwise of a
non-trivial zero.  Thus one uses an invariant $\mathcal{I}(q_1,\dots,q_r)$
of the system (as constructed by Schmidt \cite{sch}) which is a
function of the various coefficients, and one defines a minimal model to
be one in which all the forms $q_i$ are defined over $\Z_p$, and for
which the $p$-adic valuation $|\mathcal{I}(q_1,\ldots,q_r)|_p$ is
maximal.  It can happen that $\mathcal{I}(q_1,\ldots,q_r)=0$, but it is
possible to avoid consideration of such systems.

In order to get a feel for what a minimal model might look like,
observe that if one takes the transforms (\ref{xmap}) and (\ref{qmap})
to be $M=pI_n$ and $P=p^{-2}I_r$ respectively, then we return to the
original system.  Thus we may think of transforms in which 
\[|\det(P)|_p=|\det(M)|_p^{-2r/n}\]
as ``neutral''.  However if there is a
pair of transforms producing an integral system, but for which
$|\det(P)|_p>|\det(M)|_p^{-2r/n}$, then we may regard this as having
removed at least one factor $p$ from the system.  In fact the
condition for a minimal model is precisely that one should have
\begin{equation}\label{mincond}
|\det(P)|_p\le|\det(M)|_p^{-2r/n}
\end{equation}
for any transforms that produce another integral system.  We shall say that
$q_1,\ldots,q_r$ is ``minimized'' if it meets this condition.

Under the assumption that $n>4r$ we can draw certain conclusions about
the system $Q_1,\ldots,Q_r$  over $\F_p$.  For example, if there were
any form, $Q_1$, say, such that $Q_1(0,0,x_3,x_4,\ldots,x_n)$ vanishes
identically, then the transforms $M={\rm Diag}(p,p,1,1,\ldots,1)$ and
$P={\rm Diag}(p^{-1},1,1,\ldots,1)$ would map $q_1,\ldots,q_r$ to
another integral system .  However these would violate the condition
(\ref{mincond}). Thus no form $Q_i$, or more generally no form in the
linear system generated by $Q_1,\ldots,Q_r$, can be annihilated by
setting two variables to zero. In the same way one can show that one
cannot annihilate any $k$ of the forms by setting $2k$ variables to
zero. If $Q_1,\ldots,Q_r$ satisfy these conditions we will say that
the system is ``$\F_p$-minimized''.  Thus if $q_1,\dots,q_r$ is minimized,
then $Q_1,\ldots,Q_r$ is $\F_p$-minimized.  However the converse
is not true in general. For example, when $r=1$, $n=5$ and $p=3$, 
the form $q_1(x_1,\ldots,x_5)=x_1^2+x_2^2+x_3x_4+9x_5^2$ is not 
minimized, since we can take $M$ as ${\rm Diag}(1,1,1,1,1/3)$ and $P$ as the
$1\times 1$ identity matrix.  Then $|\det(P)|_3=1$, while
$|\det(M)|_3^{-2r/n}=3^{-2/5}$, contravening the condition
(\ref{mincond}). On the other hand the reduction to $\F_3$ is
$Q_1(x_1,\ldots,x_5)=x_1^2+x_2^2+x_3x_4$, which is $\F_3$ minimized.

We can illustrate the use of a minimal model by looking at the case
$r=1$, with $n\ge 5$.  For a minimal model, $Q_1$ cannot be
annihilated by setting two variables to zero.  We proceed to
make a linear change of variables so as to represent $Q_1$ using as
few variables as possible, and put
$Q_1(x_1,\ldots,x_n)=Q^*(y_1,\ldots,y_m)$ accordingly.  Then we will
have $m\ge 3$, by the minimality condition.  The
Chevalley--Warning Theorem now produces a non-trivial zero of $Q^*$.  Such
a zero must be non-singular, since otherwise the form $Q^*$ would be
degenerate, contrary to hypothesis.  This results in a non-singular
zero of $Q_1$, to which Hensel's Lemma may be applied, completing the
proof. The reader may care to note that this approach allows all the
fields $\Q_p$, including the case $p=2$, to be handled uniformly.

A similar argument handles the case $r=2$, for $n\ge 9$ (Demyanov
\cite{Dem} and Birch, Lewis and Murphy \cite{BLM}). In particular, if
$n\ge 9$ and $<Q_1,Q_2>$ is \linebreak $\F_p$-minimized, then $Q_1$ and $Q_2$
always have a non-singular common zero.  When $r=3$ and $n\ge 13$
there appear to be numerous special cases to consider.  The work of
Birch and Lewis \cite{BL} and Schuur \cite{Schuur} proves
similarly that if $<Q_1,Q_2,Q_3>$ is $\F_p$-minimized, then there is
a non-singular common zero, provided that $p\ge 11$.  
However the approach is doomed to
fail in general, as the following example shows.  We take $p=2$ and
examine the forms
\[Q_1(x_1,\ldots,x_{13})=x_1x_2+x_3^2+x_3x_4+x_4^2,\]
\[Q_2(x_1,\ldots,x_{13})=x_5x_6+x_7^2+x_7x_8+x_8^2,\]
\[Q_3(x_1,\ldots,x_{13})=x_1^2+x_1x_2+x_2^2+x_5x_7+x_6x_8+x_7^2+x_8^2\]
over $\F_2$. We claim that any common zero (over $\F_2$) is a singular
zero for $Q_1$, and hence is singular for the whole system.  To verify this
one easily checks that a non-singular zero of $Q_1$ has
$x_1^2+x_1x_2+x_2^2=1$, and that $x_5x_7+x_6x_8+x_7^2+x_8^2=0$ for any
zero of $Q_2$. This is enough to show that $Q_3=1$ at any point which
is both a nonsingular zero of $Q_1$ and a zero of $Q_2$.  The claim
then follows.  One can also verify that the system is
$\F_2$-minimized, which requires a case by case analysis.  We give a
single example, showing that 
\[Q_1+Q_3=x_1^2+x_2^2+x_3^2+x_3x_4+x_4^2+x_5x_7+x_6x_8+x_7^2+x_8^2\]
and 
\[Q_2+Q_3=x_1^2+x_1x_2+x_2^2+(x_5+x_8)(x_6+x_7)\]
cannot both vanish on a linear space $L\le \F_2^8$ of dimension 4. It
will be convenient to work with a basis $\b{e}_1,\ldots,\b{e}_8$ of
$\F_2^8$, corresponding to the variables $x_1\ldots,x_8$. If we set
$\b{e}_5'=\b{e}_5+\b{e}_8$ then on the space
\[V:=<\b{e}_1,\b{e}_2,\b{e}_3,\b{e}_4,\b{e}_5',\b{e}_6,\b{e}_7>\]
the forms become
\[x_1^2+x_2^2+x_3^2+x_3x_4+x_4^2+x_5x_6+x_5x_7+x_5^2+x_7^2\]
and
\[x_1^2+x_1x_2+x_2^2,\]
both of which must vanish on $L':=V\cap L$, which will have dimension
at least 3.  The second form vanishes only when $x_1=x_2=0$.  Hence
$L'$ must be contained in
$V':=<\b{e}_3,\b{e}_4,\b{e}_5',\b{e}_6,\b{e}_7>$. On this latter
space the first form $Q_1+Q_3$ reduces to
\[Q:=x_3^2+x_3x_4+x_4^2+x_5x_6+x_5x_7+x_5^2+x_7^2.\]
However $Q$ is non-degenerate on $V'$, and hence cannot vanish on a
subspace of dimension 3.   This contradiction shows that $Q_1+Q_3$ 
and $Q_2+Q_3$ cannot both vanish on $L$.

We therefore see that this particular line of attack cannot prove that
$\beta(r;\Q_p)=4r$ for all $r$ and $p$.  However one might consider
working modulo $p^2$ or with higher powers, instead of reducing to $\F_p$.
If one works only over $\F_p$ the following result 
seems the most that one can hope for.
\begin{theorem} (Heath-Brown \cite{comp}.)
For all $r\in\N$ one has $\beta(r;\Q_p)=4r$ if $p\ge (2r)^r$.  Indeed
if $K$ is any finite extension of $\Q_p$ with residue field $F_K$, then
$\beta(r;K)=4r$ if $\# F_K\ge (2r)^r$. Moreover an $F_K$-minimized 
system $Q_1,\ldots,Q_r$ has a
non-singular common zero provided that $\# F_K\ge (2r)^r$.
\end{theorem}

One should recall that the Ax--Kochen Theorem \cite{AK} yields
$\beta(r;Q_p)=4r$ for $p\ge p(r)$, so one might view the result above
as merely giving an explicit value for $p(r)$. However when one looks
at extensions $K$ of $\Q_p$ there is a more important difference. The
Ax--Kochen result implies that $\beta(r;K)=4r$ if the characteristic of $F_K$
is at least some value $p(r;[K:\Q_p])$.  In contrast Theorem 1
has a condition only on the size of $F_K$.  Thus it is conceivable
that the Ax--Kochen result never applies when the characteristic of 
$F_K$ is $2$, for example.

The overall plan for the proof of Theorem 1 is to give a lower bound
for the overall number of zeros of $Q_1,\ldots,Q_r$ over $F_K$, and to
compare this with an upper bound for the number of singular zeros. It
turns out that to count common zeros it suffices to count zeros of
each individual linear combination $Q_a:=a_1Q_1+\ldots+a_rQ_r$. The
number of zeros of $Q_a$ in $F_K^n$ is approximately $(\# F_K)^{n-1}$, and
the discrepancy depends (in part) on the rank of $Q_a$.  It therefore
turns out that the key step in the proof is to give a good upper bound
for the number of vectors $a$ in $F_K^r$ for which $Q_a$ has a given
rank. This step uses the minimality conditions.

An interesting corollary to Theorem 1 is provided by the following
result of Leep \cite{app}.
\begin{theorem}(Leep.)
Let $p$ be a prime and let $L=\Q_p(T_1,\ldots,T_k)$.  Then $\beta(1;L)=2^{2+k}$.
\end{theorem}
Thus the $u$-invariant of the function field $L$ is $2^{2+k}$.  Before
this result there had been much work on the case $k=1$, culminating in
a successful treatment for all primes $p\not=2$, by Parimala and
Suresh \cite{PS}. Nothing however was known for $k\ge 2$.
Now one can even handle pairs of forms, showing that
\[\beta(2;L)=2^{3+k}.\]

One striking feature of Leep's result is that, in contrast to
Theorem~1, there is no restriction on the size of $p$.  It is
interesting to see how this comes about.  Suppose a quadratic form
$q(x_1,\ldots,x_n)\in L[x_1,\ldots,x_n]$ is given. We aim to locate a zero
of $q$ in which $x_1,\ldots,x_n$ are polynomials in $T_1,\ldots,T_k$ of
degree at most $d$ say, by finding suitable values (in $\Q_p$) for the 
various coefficients $c_1,\ldots,c_N$ say. The conditions these $c_i$
have to satisfy form a system  of a large number ($R$ say) of
quadratic forms.  Here $N$ and $R$ will depend on $d$, but if
$n>2^{2+k}$ we will have $N>4R$ for large enough $d$.  Thus, by
Theorem~1, one can find suitable coefficients $c_i$ provided that
$p\ge (2R)^R$. The trick now is to use an extension
$L^*=K(T_1,\ldots,T_k)$ of $L$ obtained by taking $K$ to be an
extension of $\Q_p$ of odd degree having $\# F_K\ge (2R)^R$.
Everything now works as before, with the
values of $N$ and $R$ unaffected by this change.  The result is that
we obtain a non-trivial solution $q(x_1,\ldots,x_n)=0$ in which 
the $x_i$ are in $L^*$. Finally we appeal to a result of Springer \cite{spr},
which shows that a quadratic form over a field $L$ of characteristic
different from 2 has a non-trivial zero provided that there is a zero over
some odd degree extension of $L$. This completes the
proof. (Incidentally, although Springer's statement required
characteristic different from 2,  David Leep points out that one may 
prove the result without this restriction by essentially the same method.)
\bigskip

We turn now to the second main line of attack on $\beta(r;\Q_p)$.
This will provide upper bounds for $\beta(r;\Q_p)$ which are expected
in many cases to be sub-optimal.  However the method has the
advantage of producing results for every prime $p$. The procedure uses
induction on $r$, and originates from work of Leep \cite{leep} in
1984. One can see that the first line of attack runs into difficulties
when the field $\F_p$ is small --- it leaves too little room for
man\oe uvre.  The second strategy works purely over $\Q_p$ and so
encounters no such problems.

The basic idea is as follows.  Suppose we are given forms
$q_1,\ldots,q_r$ over $\Q_p$. If one can find a linear space $L$ in
$\Proj(\Q_p)^{n-1}$, with projective dimension $\beta(k;\Q_p)$, such that
$q_1,\dots,q_{r-k}$ all vanish identically on $L$, then the remaining
$k$ forms $q_{r-k+1},\ldots,q_r$ will have a zero in $L$, by
definition of $\beta(k;\Q_p)$.  Thus the focus of this technique is on
the number $\beta(r;K,m)$, defined as the largest integer $n$ for
which there exist quadratic forms
$q_1(x_1,\ldots,x_n),\ldots,q_r(x_1,\ldots,x_n)$ over $K$ such that
there is no $K$-linear space of projective dimension $m$ on which the
forms vanish identically.  The argument above shows now that
\begin{equation}\label{ind1}
\beta(r;K)\le\beta(r-k;K,\beta(k;K))
\end{equation}
for $k<r$, for any field $K$.

One may estimate $\beta(r;K,m)$ via induction on $m$. Suppose our
forms vanish on a projective linear space $L$ of dimension $m-1$,
spanned by $\b{e}_0,\ldots,\b{e}_{m-1}$ say. We wish to find an
additional vector $\b{e}_m=\b{e}$ to add to this basis.\linebreak  Let $L^*$ be
a complementary linear space for $L$ in $\Proj(K)^{n-1}$, so that
$\dim(L^*)=n-m-1$. We will require $[\b{e}]$ to belong to $L^*$, which
will ensure that $\b{e}_0,\ldots,\b{e}_{m-1},\b{e}$ are linearly
independent. In order for our forms to vanish on the span of the
extended set $\b{e}_0,\ldots,\b{e}_{m-1},\b{e}$ it suffices that
\[q_i(\b{e}_j,\b{e})=0,\;\;\;(1\le i\le r,\; 0\le j\le m-1)\]
and
\[q_i(\b{e})=0,\;\;\;(1\le i\le r)\]
where $q_i(\x,\b{y})$ is the bilinear form associated to $q_i$.
The first set of conditions restricts $\b{e}$ to a subspace of $L^*$ of
codimension at most $rm$, so that a suitable $\b{e}$ must exist,
provided that $\beta(r;K)<n-m-rm$.  It follows that our basis can be
extended whenever $n>(r+1)m+\beta(r;K)$, yielding the inductive
inequality 
\[\beta(r;K,m)\le\max\left\{\beta(r;K,m-1)\,,\,(r+1)m+\beta(r;K)\right\}.\]
We therefore deduce that
\begin{equation}\label{ind2}
\beta(r;K,m)\le(r+1)m+\beta(r;K)
\end{equation}
for all $m$.

One may combine this with (\ref{ind1}) to obtain
\[\beta(r;\Q_p)\le\beta(r-2;\Q_p,\beta(2;\Q_p))=\beta(r-2;\Q_p,8)
\le 8(r-1)+\beta(r-2;\Q_p).\]
Starting from $\beta(1;\Q_p)=4$ and $\beta(2;\Q_p)=8$ one then finds
that
\begin{equation}\label{mr}
\beta(r;\Q_p)\le\left\{\begin{array}{cc} 2r^2, & r \mbox{ even},\\
2r^2+2, & r \mbox{ odd},
\end{array}\right.
\end{equation}
(Martin \cite{mar}, improving slightly on the original result of Leep).
In particular one has
\[\beta(3;\Q_p)\le 20,\]
for all primes $p$.

One may ask whether one can improve on the bound (\ref{ind2}).  In the
case $r=1$ the estimate (\ref{ind2}) becomes 
$\beta(1;\Q_p,m)\le 2m+4$, and indeed this is
best possible. However for $r=2$ one has only $\beta(2;\Q_p,m)\le
3m+8$, and here one can do better by an argument due to Dietmann
\cite{Diet} (improved slightly by Heath-Brown \cite{hbquartics}).  The
method is based on the following theorem of Amer \cite[Satz 8,
p.29]{amer} in an unpublished thesis.
\begin{theorem}(Amer, 1976) For any field $K$ of characteristic
  $\chi_K\not=2$ one has $\beta(2;K,m)\le \beta(1; K(T),m)$ for every
  integer $m\ge 0$.
\end{theorem}
The special case $m=0$ is given by Brumer \cite{Brum}. In an
unpublished manuscript Leep shows that the result holds even when
$\chi_K=2$. 

In view of (\ref{ind2}) one has $\beta(1; K(T),m)\le
2m+\beta(1;K(T))$, so that
\[\beta(2;K,m)\le 2m+\beta(1;K(T)).\]
We take $K=\Q_p$ and use the case $k=1$ of Theorem 2, which produces 
$\beta(1;\Q_p(T))=8$.  We therefore conclude that
\begin{equation}\label{d2}
\beta(2;\Q_p,m)\le 2m+8
\end{equation}
for all $m\ge 0$, which is easily shown to be best possible.
Unfortunately it seems that we can get results of this quality only
for the cases $r=1$ and $r=2$. We therefore ask:
\bigskip

{\bf Open Question} {\em Is it true that} $\beta(3;\Q_p,m)=2m+O(1)$
{\em uniformly for all} $m\ge 1$ {\em and all primes} $p$ ?
\bigskip

Even the situation over $\F_p$ is unclear. One can use Amer's theorem
to show that $\beta(2;\F_p,m)=2m+4$ for primes $p\ge 3$, and the
result of Leep noted above similarly handles $p=2$.
However it appears to be unknown whether or not $\beta(3;\F_p,m)=2m+O(1)$.

We can use (\ref{d2}) to advantage in our previous argument.
From (\ref{ind1}) we have $\beta(r;\Q_p)\le
\beta(2;\Q_p,\beta(r-2;\Q_p))$, so that (\ref{d2}) yields
\[\beta(r;\Q_p)\le 2\beta(r-2;\Q_p)+8.\]
In particular
\begin{equation}\label{d3}
\beta(3;\Q_p)\le 2\beta(1;\Q_p)+8=8+8=16,
\end{equation}
\[\beta(4;\Q_p)\le 2\beta(2;\Q_p)+8=24,\]
\[\beta(5;\Q_p)\le 2\beta(3;\Q_p)+8\le 40\]
using (\ref{d3}), and
\[\beta(6;\Q_p)\le 2\beta(4;\Q_p)+8\le 56.\]
One may then use (\ref{ind1}) with $k=1$ and (\ref{ind2}) with 
$m=\beta(1;\Q_p)=4$ to show that
\[\beta(7;\Q_p)\le \beta(6;\Q_p)+28\le 84.\]
From this point on the most efficient procedure is to use
(\ref{ind1}) with $k=2$ and (\ref{ind2}) with $m=\beta(2;\Q_p)=8$,
deducing that
\[\beta(r;\Q_p)\le\left\{\begin{array}{cc} 2r^2-16, & r \mbox{ even }
    \ge 6,\\
2r^2-14, & r \mbox{ odd } \ge 7,
\end{array}\right.\]
which improves on Martin's result (\ref{mr}) by 16. Thus the overall saving
is not large, but is not insignificant for $r=3$, for example.

In the simplest open case $r=3$ our state of knowledge is therefore that
\[12\le\beta(3;\Q_p)\le 16.\]
It is perhaps of interest to review the somewhat roundabout route to
the upper bound here, since it combines results from both the lines of
attack described here. The steps could be summarized as follows.
\begin{enumerate}
\item Theorem 1 handles systems of $r$ forms in at least $4r+1$
  variables, over an extension $K$ of $\Q_p$, when $\# F_K\ge (2r)^r$.
\item Leep's argument for the proof of Theorem 2 shows that 
if 
\[q(x_1,\ldots,x_9)\in\Q_p(T)[x_1,\ldots,x_9]\]
then one can find a zero
over an extension $K(T)$, provided that $\# F_K$ is large enough.
\item On choosing a suitable odd degree extension, Springer's theorem
  produces a zero of $q(x_1,\ldots,x_9)$ over $\Q_p(T)$.
\item For suitable $n$ we may then find a large linear space of 
solutions for a form
$q(x_1,\ldots,x_n)$ over $\Q_p(T)$.
\item Amer's theorem then shows that a pair of forms
  $q_i(x_1,\ldots,x_n)$ (for $i=1,2$) over $\Q_p$ also has a large
  linear space of solutions.
\item Combined with the estimate (\ref{ind1}) this produces our bound
  for $\beta(3;\Q_p)$.
\end{enumerate}

It certainly seems surprising that the proof goes via systems of large
numbers of forms in steps 1 and 2.  It would be interesting to know whether
an argument based on minimal models could give a direct proof of a
bound weaker than $\beta(3;\Q_p)=12$, but valid for all $p$.

\bigskip
\bigskip

Mathematical Institute,

24--29, St. Giles',

Oxford

OX1 3LB

UK
\bigskip

{\tt rhb@maths.ox.ac.uk}


\begin{thebibliography}{9}
\bibitem{amer}M. Amer, {\em Quadratische formen \"{u}ber 
funktionenk\"{o}rpern}, (Thesis, Mainz, 1976).
\bibitem{Artin}E. Artin, {\em The collected papers of Emil Artin}, 
(Addison--Wesley, London, 1965).
\bibitem{AK}J. Ax and S. Kochen, 
Diophantine problems over local fields. I,
{\em Amer. J. Math.}, 87 (1965), 605--630.
\bibitem{birch} B.J. Birch, Forms in many variables,
{\em Proc. Roy. Soc. Ser. A} 265 (1961/1962), 245--263.
\bibitem{BLM}B.J. Birch, D.J. Lewis and T.G. Murphy, 
Simultaneous quadratic forms, {\em Amer. J. Math.,} 84 (1962), 110--115. 
\bibitem{BL}B.J. Birch and D.J. Lewis, Systems of three quadratic
forms, {\em Acta Arith.,} 10 (1964/1965), 423--442.
\bibitem{Brum}A. Brumer, Remarques sur les couples de formes 
quadratiques, {\em C. R. Acad. Sci. Paris S\'er. A-B,} 286 (1978),
no. 16, A679--A681.
\bibitem{Dem}V.B. Demyanov, Pairs of quadratic forms over a complete 
field with discrete norm with a finite field of residue classes,
{\em Izv. Akad. Nauk SSSR. Ser. Mat.,} 20 (1956), 307--324.
\bibitem{Diet}R. Dietmann, Linear spaces on the intersection of 
two quadratic hypersurfaces, and systems of $p$-adic quadratic 
forms, {\em Monatsh. Math.}  146  (2005), 175--178.
\bibitem{Has}H. Hasse, Darstellbarkeit von Zahlen durch quadratische 
Formen in einem beliebigen algebraischen
Zahlk\"{o}rper, {\em J. Reine Angew. Math.,} 153 (1924), 11--130.
\bibitem{hbquartics}D.R. Heath-Brown, Zeros of $p$-adic forms, {\em 
Proc. Lond. Math. Soc. (3)} 100 (2010), 560--584.
\bibitem{comp}D.R. Heath-Brown, Zeros of systems of 
${\mathfrak p}$-adic quadratic forms,
{\em Compositio Math.,} 146 (2010), 271-287. 
\bibitem{leep}D.B. Leep, Systems of quadratic forms,
{\em J. Reine Angew. Math.}, 350 (1984), 109--116.
\bibitem{app}D.B. Leep, The $u$-invariant of $p$-adic function fields,
{\em J. Reine Angew. Math.}, to appear.
\bibitem{mar}G. Martin, Solubility of systems of quadratic forms,
{\em Bull. London Math. Soc.,} 29 (1997), 385--388. 
\bibitem{PS}R. Parimala and V. Suresh, 
The $u$-invariant of the function fields of $p$-adic curves,
{\em Ann. of Math. (2)}, 172 (2010), 1391--1405. 
\bibitem{sch}W.M. Schmidt, Simultaneous $p$-adic zeros of quadratic forms,
{\em Monatsh. Math.}, 90 (1980), 45--65. 
\bibitem{Schuur}S.E. Schuur, On systems of three quadratic forms,
{\em Acta Arith.,} 36 (1980), 315--322. 
\bibitem{spr}T.A. Springer, Sur les formes quadratiques d'indice 
z\'{e}ro, {\em C. R. Acad. Sci. Paris}, 234 (1952), 1517--1519. 
\bibitem{Ter}G. Terjanian,  Un contre-exemple \`a une conjecture 
d'Artin, {\em  C. R. Acad. Sci. Paris S\'er.} A-B,  262  (1966), A612.
\end{thebibliography}
\end{document}